\documentclass[twocolumn,prb,showpacs]{revtex4-1}
\usepackage{graphicx}
\usepackage{epstopdf}
\usepackage{dcolumn}
\usepackage{amsmath}
\usepackage{amsfonts}
\usepackage{multirow} 
\usepackage[active]{srcltx}
\usepackage{hyperref}
\usepackage{yfonts}

\begin{document}

 \title{The Mayan Long Count Calendar} 

\author{Thomas Chanier$^{*1}$}
\affiliation{$^1$ Universit\'e de Bretagne Occidentale, 6 avenue Victor le Gorgeu, F-29285 Brest Cedex, France}


\begin{abstract}
The Mayan Codices, bark-paper books from the Late Postclassic period
(1300 to 1521 CE) contain many astronomical tables correlated to ritual cycles, evidence of the achievement of  Mayan naked-eye astronomy and mathematics in connection to religion. In this study, a calendar supernumber is calculated by computing the least common multiple of 8 canonical astronomical periods. The three major calendar cycles, the Calendar Round, the Kawil and the Long Count Calendar are shown to derive from this supernumber. The 360-day Tun, the 365-day civil year Haab' and the 3276-day Kawil-direction-color cycle are determined from the prime factorization of the 8 canonical astronomical input parameters. The 260-day religious year Tzolk'in and the Long Count Periods (the 360-day Tun, the 7200-day Katun and the 144000-day Baktun) result from arithmetical calculations on the calendar supernumber. My calculations explain certain aspect of the Mayan Calendar notably the existence of the Maya Epoch, a cycle corresponding to 5 Maya Eras of 13 Baktun. Modular arithmetic considerations on the calendar supernumber give the position of the Calendar Rounds at the Mayan origin of time, the Long Count Calendar date 0.0.0.0.0 4 Ahau 8 Cumku. Various long count numbers identified on Mayan Codices and monuments are explained by my approach. In particular, the results provide the meaning of the Xultun numbers, four enigmatic long count numbers deciphered in 2012 by Saturno {\it et al.} on the inner walls of a small masonry-vaulted structure in the extensive Mayan ruins of Xultun, Guatemala.\cite{saturno} The results show a connection between the religious sites of Xultun and Chich\'en Itz\'a, Mexico. This is the first study linking unambiguously Mayan astronomy, mathematics and religion, providing a unified description of the Mayan Calendar.
\end{abstract}

\maketitle


Mayan priests-astronomers were known for their astronomical and mathematical proficiency, as exemplified in the Dresden codex, a bark-paper book of the XI or XII century CE containing many astronomical tables correlated to ritual cycles. However, due to the zealous role of the Inquisition during the XVI century CE Spanish conquest of Mexico, number of these Codices were burnt, leaving us with few information on Pre-Columbian Mayan culture. Thanks to the work of Mayan archeologists and epigraphists since the early XX century, the few codices left, along with numerous inscriptions on monuments, were deciphered, underlying the importance of the concept of time in Maya civilisation. This is reflected by the three major Mayan Calendars, reminiscent of the Mayan cyclical conception of time: the Calendar Round, the Long Count Calendar and the Kawil-direction-color cycle.

\begin{figure}[h!]
\begin{center}
\includegraphics[width = 65mm]{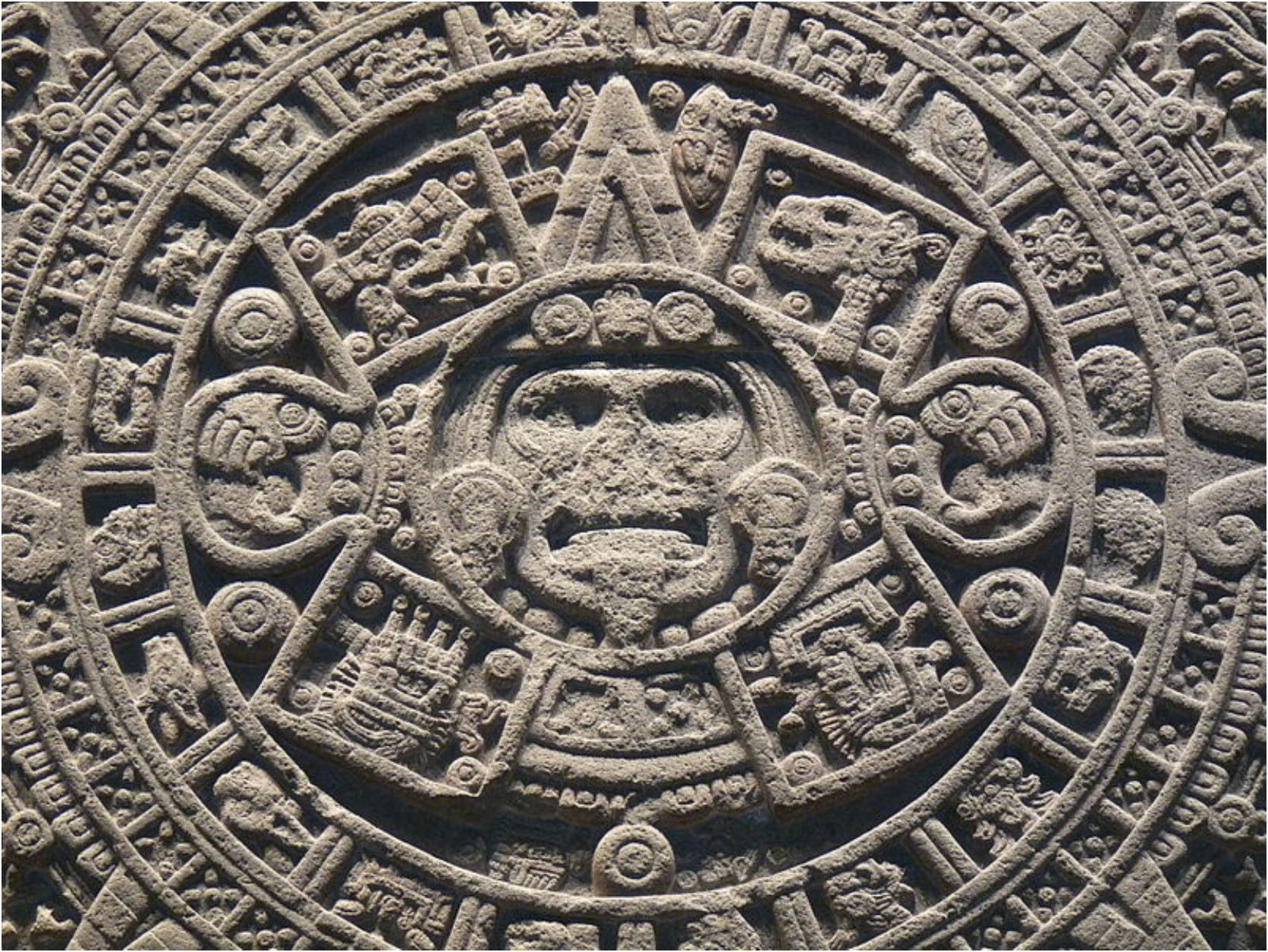}
\caption{Mayan/Aztec Calendar stone representing the Five Suns, discovered in 1790 at El Z\'ocalo, Mexico City, Mexico.}\label{YYY}
\end{center}
\end{figure}

The Calendar Round (CR) represents a day in a non-repeating 18980-day cycle, a period of roughly 52 years, the combination of the 365-day civil year Haab' and the 260-day religious year Tzolk'in. The Tzolk'in comprises 13 months (numerated from 1 to 13) containing 20 named days (Imix, Ik, Akbal, Kan, Chicchan, Cimi, Manik, Lamat, Muluc, Oc, Chuen, Eb, Ben, Ix, Men, Cib, Caban, Etznab, Cauac, and Ahau). This forms a list of 260 ordered Tzolk'in dates from 1 Imix, 2 Ik, ... to 13 Ahau.\cite{AveniTz} The Haab' comprises 18 named months (Pop, Uo, Zip, Zotz, Tzec, Xul, Yaxkin, Mol, Chen, Yax, Zac, Ceh, Mac, Kankin, Muan, Pax, Kayab, and Cumku) with 20 days each (Winal) plus 1 extra month (Uayeb) with 5 nameless days. This forms a list of 365 ordered Haab' dates from 0 Pop, 1 Pop, ... to 4 Uayeb.\cite{AveniHa} The Tzolk'in and the Haab' coincide every 73 Tzolk'in or 52 Haab' or a Calendar Round such as the least common multiple (LCM) of 260 and 365: $73\times260=52\times365=18980$ days. In the Calendar Round, a date is represented by $\alpha$X$\beta$Y with the religious month $1\leq\alpha\leq13$, X one of the 20 religious days, the civil day $0\leq\beta\leq19$, and Y one of the 18 civil months, $0\leq\beta\leq4$ for the Uayeb. 
To reckon time in a period longer than 52 years, the Maya used the Long Count Calendar (LCC), describing a date $D$ in a 1872000-day Maya Era of 13 Baktun, a period of roughly 5125 years, counting the number of day elapsed from the Mayan origin of time, the mythical date of creation 0.0.0.0.0 4 Ahau 8 Cumku, corresponding to the Gregorian Calendar date 11 August 3114 BC according to the Goodman-Martinez-Thompson correlation. The XXI century saw the passage of the new Era on 21 December 2012, a date related to several apocaliptic or world renewal New Age theories, relayed by mass-media. Whereas Mayan mathematics are based on a vigesimal basis, the LCC is a 18-20-mixed radix basis representation of time: a date $D$ is defined by a set of 5 numbers ($C_4.C_3.C_2.C_1.C_0$) such as $D\equiv \mathrm{mod}(D,13 \times 144000) = C_0+C_1\times20+C_2\times360+C_3\times7200+C_4\times144000$ where $C_4$ is the number of Baktun (144000 days) in the 13 Baktun Era ($0\leq C_4<13$), $C_3$ the number of Katun (7200 days) in the current Baktun ($0\leq C_3<20$), $C_2$ the number of Tun (360 days) in the current Katun ($0\leq C_2<18$), $C_1$ the number of Winal (20 days) in the current Tun ($0\leq C_1<20$) and $C_0$ the number of Kin (days) in the current Winal ($0\leq C_0<20$). The Kawil-direction-color cycle or 4-Kawil is a 3276-day cycle, the combination of the 4 directions-colors and the 819-day Kawil.\cite{Berlin} Table \ref{CC} gives the different  calendar cycles with their prime factorizations.

\begin{table}[h!]
\begin{center}
\begin{tabular}{l l l}
Cycle & $P$ [day] & Prime factorization\\\hline
Haab' & 365 & 5 $\times$ 73\\
Tzolk'in & 260 & 2$^2$ $\times$ 5 $\times$ 13\\\hline
Winal & 20 & 2$^2$ $\times$ 5\\
Tun & 360 & 2$^3$ $\times$ 3$^2$ $\times$ 5\\
Katun & 7200 & 2$^5$ $\times$ 3$^2$ $\times$ 5\\
Baktun & 144000 & 2$^7$ $\times$ 3$^2$ $\times$ 5\\\hline
Kawil & 819 & 3$^2$ $\times$ 7 $\times$ 13\\
4-Kawil & 3276 & 2$^2$ $\times$ $3^2$ $\times$ 7 $\times$ 13\\\hline
\end{tabular}
\caption{Calendar cycles and their prime factorizations.}\label{CC}
\end{center}
\end{table}

The origin of the Long Count Periods is unknown. A common assumption is the desire of the calendar keeper to maintain the Tun in close agreement with the tropical/solar year of approximately 365.24 days.\cite{aveni1} There is no consensus concerning the origin of the Tzolk'in, which has been associated with various astronomical cycles. 3 Tzolk'in correspond to Mars synodic period, 16 Tzolk'in equal 11 of Saturn synodic periods (+2 days), and 23 Tzolk'in are equivalent to 15 Jupiter synodic periods (-5 days). \cite{justeson} It has been tentatively connected to the eclipse half-year (173.31 days) because 2 Tzolk'in are very close to 3 eclipse half-years.\cite{aveni2} Finally, it has been noted that the Tzolk'in approximates the length of time Venus is visible as a morning or evening star. \cite{bricker69} The Kawil cycle has been attributed to the observation of Jupiter and Saturn\cite{lounsbury1,justeson1} because 19 (6) Kawil correspond to 39 (13) Jupiter (Saturn) synodic period. However, these interpretations fail to link the Tzolk'in and the Kawil to the Long Count Periods. The reason why the initial state of the Calendar Round at the LCC origin of time 0.0.0.0.0 is 4 Ahau 8 Cumku remains unexplained up to now.

In 2012, four LCC numbers, the Xultun numbers (Table \ref{Xultun}), have been discovered on the walls of a small painted room in the Mayan ruins of Xultun, dating from the early IX century CE.\cite{saturno} These numbers have a potential astronomical meaning. Indeed, $\mathcal{X}_0$ = LCM(260,360,365) is a whole multiple of the Tzolk'in, Haab', Tun, Venus and Mars synodic periods: 341640 = 1314 $\times$ 260 = 936 $\times$ 365 = 949 $\times$ 360 = 585 $\times$ 584 = 438 $\times$ 780, $\mathcal{X}_1=365\times3276$ is the commensuration of the Haab' and the 4-Kawil cycle. The greatest common divisor of the $\mathcal{X}_i$'s is 56940 = 3 CR corresponding to the commensuration of the Haab' and Mars synodic period of 780 days. However, the meaning of $\mathcal{X}_2$ and $\mathcal{X}_3$, related to $\mathcal{X}_0$ by the equation $\mathcal{X}_3=\mathcal{X}_2+2\mathcal{X}_0$ is unknown.

\begin{table}[h!]
\begin{center}
\begin{tabular}{l l l c l l}
$\mathcal{X}_i$ & LCC & $D$ [day] & $\mathcal{X}_i/56940$ & \\\hline
$\mathcal{X}_0$ &	2.7.9.0.0 & 341640 & 6\\
$\mathcal{X}_1$ &	8.6.1.9.0 & 1195740 & 21\\
$\mathcal{X}_2$ &	12.5.3.3.0 & 1765140 & 31\\
$\mathcal{X}_3$ &	17.0.1.3.0 & 2448420 & 43\\
\hline
\end{tabular}
\caption{Xultun numbers $\mathcal{X}_i$.\cite{saturno} 56940 = LCM(365,780) = 2$^2$ $\times$ 3 $\times$ 5 $\times$ 13 $\times$ 73 is their largest common divisor and $\mathcal{X}_3=\mathcal{X}_2+2\mathcal{X}_0$.}\label{Xultun}
\end{center}
\end{table}

\begin{table}[h!]
\begin{center}
\begin{tabular}{l l c l}
Planet &$i$ &  $P_i$ [day] & Prime factorization\\\hline
Mercury &1 & 116 & 2$^2$ $\times$ 29\\
Venus &2 & 584 & 2$^3$ $\times$ 73\\
Earth &3 & 365 & 5 $\times$ 73\\
Mars &4 & 780 & 2$^2$ $\times$ 3 $\times$ 5 $\times$ 13\\
Jupiter &5 & 399 & 3 $\times$ 7 $\times$ 19\\
Saturn &6 & 378 & 2 $\times$ 3$^3$ $\times$ 7\\\hline
Lunar &7 & 177 & 3 $\times$ 59\\
senesters &8 &  178 & 2 $\times$ 89\\\hline

\end{tabular}
\caption{Planet canonical cycles \cite{saturno,NASA,bricker} and their prime factorizations.}\label{SP}
\end{center}
\end{table}

Mayan astronomers-priests, known for their astronomical proficiency, may have observed with a naked eye the periodic movements of the five planets visible in the night sky:  the moon, Mercury, Venus, Earth (solar year), Mars, Jupiter, and Saturn. Their respective canonical synodic periods are given in Table \ref{SP}. Evidences have been found in different Mayan Codices for Mercury, Venus, and Mars, but it is unclear whether they tracked the movements of Jupiter and Saturn.\cite{bricker1} In particular, on page 24 of the Dresden codex is written the so-called Long Round number noted 9.9.16.0.0 or 1366560 days, a whole multiple of the Tzolk'in, the Haab', the Tun, Venus and Mars synodic periods, the Calendar Round and the Xultun number $\mathcal{X}_0$: $\cal{LR}$ = 1366560 = 5256 $\times$ 260 = 3744 $\times$ 365 = 3796 $\times$ 360 = 2340 $\times$ 584 = 1752 $\times$ 780 = 72 $\times$ 18980 = 4 $\times$ 341640. The relevant periods for the prediction of solar/lunar eclipses are the lunar semesters of 177 or 178 days (6 Moon synodic periods), which are the time intervals between subsequent eclipse warning stations present in the Eclipse Table in the Dresden Codex and the lunar tables inscribed on the Xultun walls.\cite{saturno,bricker} From their prime factorizations (Table \ref{SP}), we calculate the calendar supernumber $\mathcal{N}$ defined as the least common multiple of the $P_i$'s:

\begin{eqnarray}\label{astro}
\mathcal{N} & = & 20757814426440\\\nonumber
 & = & 2^2 \times 3^2\times 7 \times 13 \times 19 \times 29 \times 59 \times 73 \times 89\\\nonumber
 & = & 365 \times 3276 \times 2 \times 3 \times 19 \times 29 \times 59 \times 89\\\nonumber
 & = & \mathrm{LCM}(360,\ 365,\ 3276) \times 3 \times 19 \times 29 \times 59 \times 89
\end{eqnarray}

Equation \ref{astro} gives the calendar supernumber and its prime factorization. It is expressed as a function of the Haab' and the 4-Kawil. The Haab', canonical solar year, is such as the Haab' and the $P_i$'s are relatively primes (exept Venus and Mars): the \{LCM($P_i$,365)/365, $i$ = 1..8\} = \{116, 8, 1, 156, 399, 378, 177, 178\} (Table \ref{SP}). The 4-Kawil and the Haab' are relatively primes: the LCM(365,3276) = 365 $\times$ 3276 and their largest common divisor is 1. The 4-Kawil has the following properties: the \{LCM($P_i$,3276)/3276, $i$ = 1..8\} = \{29, 146, 365, 5, 19, 3, 59, 89\}. That defines the 4-Kawil. The commensuration of the 4-Kawil and the Haab' $\mathcal{X}_1$ = 365 $\times$ 3276 gives: \{LCM($P_i$,$\mathcal{X}_1$)/$\mathcal{X}_1$, $i$ = 1..8\} = \{29, 2, 1, 1, 19, 3, 59, 89\}. 360 is the integer closest to 365 such that the LCM(360,3276) = 32760 and the \{LCM($P_i$,32760)/32760, $i$ = 1..8\} = \{29, 73, 73, 1, 19, 3, 59, 89\}. The number 32760 or 4.11.0.0 has been derived from inscriptions on the Temple of the Cross in Palenque, Chiapas, Mexico.\cite{lounsburyPa} The Tun has the following properties: $\mathcal{Y}$ = LCM(360,365,3276) = 2391480 such as \{LCM($P_i$,$\mathcal{Y}$)/$\mathcal{Y}$, $i$ = 1..8\} = \{29, 1, 1, 1, 19, 3, 59, 89\}. The commensuration of the Haab', the 4-Kawil and the Baktun (400 Tun = 144000 days) gives rise to a calendar Grand Cycle $\mathcal{GC}$ = LCM(365,3276,144000) = 400 $\times$ LCM(360,365,3276) = 7 $\times$ 400 $\times$ $\mathcal{X}_0$ = 956592000. The Euclidean division of the calendar supernumber $\mathcal{N}$ by $\mathcal{GC}$ gives:

\begin{eqnarray}\label{astro2}
\mathcal{N}&=&\mathcal{GC} \times \mathcal{Q}+ \mathcal{R}\\\nonumber
\mathcal{Q} &=& 21699 \\\nonumber
\mathcal{R} &=& 724618440\\\nonumber
&=& 101 \times126\times56940\\\nonumber
 &= & 126\times\sum_{i=0}^3\mathcal{X}_i.
\end{eqnarray}

If we pose  $\mathcal{A}$ = 13 $\times$ 73 $\times$ 144000 = 400 $\times$ $\mathcal{X}_0$ = 100 $\times$ $\mathcal{LR}$ such as $\mathcal{GC}$ = 7 $\times$ $\mathcal{A}$, the Euclidean division of $\mathcal{N}$ by $\mathcal{A}$ gives:

\begin{eqnarray}\label{astro1}
\mathcal{N} &=&  \mathcal{A}  \times\mathcal{Q} +\mathcal{R}\\\nonumber
\mathcal{Q} &=& 151898\\\nonumber
\mathcal{R} &=& 41338440 = 121 \times 341640\\\nonumber &=& 121\times\mathcal{X}_0
\end{eqnarray}

The only $P_i$'s commensurate to $\mathcal{A}$ = 13 $\times$ 73 $\times$ 144000 are the Haab', Venus and Mars canonical periods (Table \ref{SP}), such as the LCM(584,365) = 37960 and LCM(780,365) = 56940 the length of the Venus and Mars Table in the Dresden codex. We have $\mathcal{A}$ = LCM(260,365,144000) = 100 $\times$ $\mathcal{LR}$ = 7200 $\times$ 18980 =  3600 $\times$ 37960 = 2400 $\times$ 56940 = 1000 $\times$ 234 $\times$ 584 = 60 $\times$ 2920 $\times$ 780. The commensuration of the Winal and 234 is the LCM(20,234) = 2340 = 9 $\times$ 260 = 20 $\times$ 117 = LCM(9,13,20). This 2340-day cycle is present in the Dresden Codex on pages D30c-D33c and has been attributed to a Venus-Mercury almanac because 2340 = 20 $\times$ 117 = 5 $\times$ 585 is an integer multiple of Mercury and Venus mean synodic periods (+1 day).\cite{bricker235} Another explanation may be of divination origin because 117 = 9 $\times$ 13. In Mesoamerican mythology, there are a set of 9 Gods called the Lords of the Night \cite{aveni,boone,lounsbury,dioses} and a set of 13 Gods called the Lords of the Day.\cite{dioses} Each day is linked with 1 of the 13 Lords of the Day and 1 of the 9 Lords of the Night in a repeating 117-day cycle. We can rewrite Equ. \ref{astro2} and \ref{astro1} as:

\begin{eqnarray}\label{astro3}
\mathcal{N} &-& 121 \times \mathcal{X}_0 = 151898\times \mathcal{A} \\\nonumber
\mathcal{N} &-& 126 \times  \sum_{i=0}^3\mathcal{X}_i = 151893\times \mathcal{A} 
\end{eqnarray}

The Long Count Periods appear in Equ. \ref{astro3}: 151898 = 338 + 360 + 7200 + 144000 and 151893 = 333 + 360 + 7200 + 144000. Adding the two equations in \ref{astro3}, we obtain:

\begin{eqnarray}\label{astro4}
5\times\mathcal{A} &=&5  \times\mathcal{X}_0 + 95 \times 126 \times 56940\\\nonumber
5\times\mathcal{A} &=&5  \times\mathcal{X}_0 + \mathrm{LCM}(\mathcal{X}_1 + \mathcal{X}_2 + \mathcal{X}_3,\ \mathcal{X}_1 + 2 \mathcal{X}_2 + \mathcal{X}_3)
\end{eqnarray}

Since $\mathcal{X}_1$ is known (the commensuration of the Haab' and the 4-Kawil), that defines  $\mathcal{X}_2$ and  $\mathcal{X}_3$. The relationship $\mathcal{X}_3=\mathcal{X}_2+2\mathcal{X}_0$ may be a mnemotechnic tool to calculate $\mathcal{Y}=\mathcal{X}_3+2\mathcal{X}_2=105\times56940$. The four Xultun numbers provides a proof that Mayan astronomers-priests determined the canonical synodic periods of the five planets visible with a naked eye: Mercury, Venus, Mars, Jupiter and Saturn. The Mayan Calendar is then constituted by a calendar Grand Cycle $\mathcal{GC}$ = 7 $\times$ $\mathcal{A}$ such as 5 $\times$ $\mathcal{A}$  = 5 $\times$ 13 $\times$ 73 $\times$ 144000 = 12000 $\times$ 56940 = 73 $\times$ $\mathcal{E}$ with $\mathcal{E}$ = 5 $\times$ 13 $\times$ 144000 is the Maya Epoch corresponding to 5 Maya Eras of 13 Baktun. This corresponds to the interpretation of the Aztec Calendar stone (Fig. \ref{YYY}): the four squares surrounding the central deity represent the four previous suns or Eras and the center deity represents the Fifth Sun or the common Era. The Tzolk'in is defined by the LCM(13,20). It is such that the LCM(260,365) is 73 $\times$ 260 = 52 $\times$ 365 = 18980 days or a Calendar Round.

A question arises at this point to know how the Maya, as early as the IX century CE, were able to comput tedious arithmetical operations on such large numbers with up to 14 digits in decimal basis. Here is a possible method. They determined the prime factorizations of the canonical astronomical periods $P_i$ (Table \ref{SP}) and listed each primes $p_i$ with their maximal order of multiplicity $\alpha_i$. They calculated the calendar supernumber $\mathcal{N}$ (the LCM of the $P_i$'s) by multiplying each $p_i$'s $\alpha_i$ time. The Euclidean division of $\mathcal{N}$ by $\mathcal{GC}$ = 7 $\times$ $\mathcal{A}$ = 7 $\times$ 400 $\times$ $\mathcal{X}_0$ (Equ. \ref{astro2}) is equivalent to a simplification of $\mathcal{N}$ by 7 $\times$ 341640 and the Euclidian division of the product of the 5 left primes (3 $\times$ 19 $\times$ 29 $\times$ 59 $\times$ 89 = 8679903) by 400. The Euclidean division of $\mathcal{N}$ by 13 $\times$ 73 $\times$ 144000 = 400 $\times$ $\mathcal{X}_0$ (Equ. \ref{astro1}) is equivalent to a simplification of $\mathcal{N}$ by 341640 and the Euclidian division of the product of the 6 left primes (3 $\times$ 7 $\times$ 19 $\times$ 29 $\times$ 59 $\times$ 89 = 60759321) by 400. It is to be noted that the prime factorization of the calendar supernumber only includes prime numbers $<100$ which facilitates the operation (there are only 25 prime numbers lower than 100). 

Modular arithmetic considerations on the calendar supernumber allows to calculate the position of the Calendar Round at the Mayan mythical date of creation. For that purpose, we first create ordered lists of the Haab' and the Tzolk'in, assigning a number on the Haab' month and the Tzolk'in day.\cite{AveniTz,AveniHa} For the Haab', the first day is 0 Pop (numbered 0) and the last day 4 Uayeb (numbered 364). For the Tzolk'in, the first day is 1 Imix (numbered 0) and the last day 13 Ahau (numbered 259). In this notation, the date of creation 4 Ahau 8 Cumku is equivalent to \{160;349\} and a date $D$ in the Calendar Round can be written as $D\equiv$ \{mod($D$ + 160,260);mod($D$ + 349,365)\}. The Calendrical Supernumber is such that: mod($\mathcal{N}$/13/73,260) = 160, mod($\mathcal{N}$/13/73,13) = 4, mod($\mathcal{N}$/13/73,20) = 0 and mod($\mathcal{N}$/13/73,73) = 49. The choice of 73 instead of 365 may be because the largest common divisor of 260 ($2^2\times5\times13$) and 365 (5 $\times$ 73) is 5. The date \{160;49\} corresponds to 4 Ahau 8 Zip, the day 0 (mod 18980), the beginning/completion of a Calendar Round. We now consider the Xultun number $\mathcal{X}_0$ = LCM(260,365,360) = 18 $\times$ 18980  = 341640 which corresponds to the time interval between two days of the same Haab', Tzolk'in, and Tun date, for example between the date origin 0.0.0.0.0 4 Ahau 8 Cumku and the date 2.7.9.0.0 4 Ahau 8 Cumku corresponding also to a completion of a 13 Tun cycle, a period of 4680 days. The completion of a Calendar Round corresponds to 18980 days elapsed such as mod(18980,4680) = 260. Starting the CR count at 4 Ahau 8 Zip, the next date in the ordered CR list such as mod($D$,4680) = 0 is the date 4 Ahau 8 Cumku, 4680 days later. A date $D$  is then expressed as \{mod($D$ + 4680 + 160,260);mod($D$ + 4680 + 49,365)\} = \{mod($D$ + 160,260);mod($D$ + 349,365)\}. The Calendar Round started therefore on a 4 Ahau 8 Zip, 4680 days earlier than the Long Count Calendar such that the starting date of the LCC is 0.0.0.0.0 4 Ahau 8 Cumku  \{160;349\}. 

\begin{figure}[h!]
\includegraphics[width = 80mm]{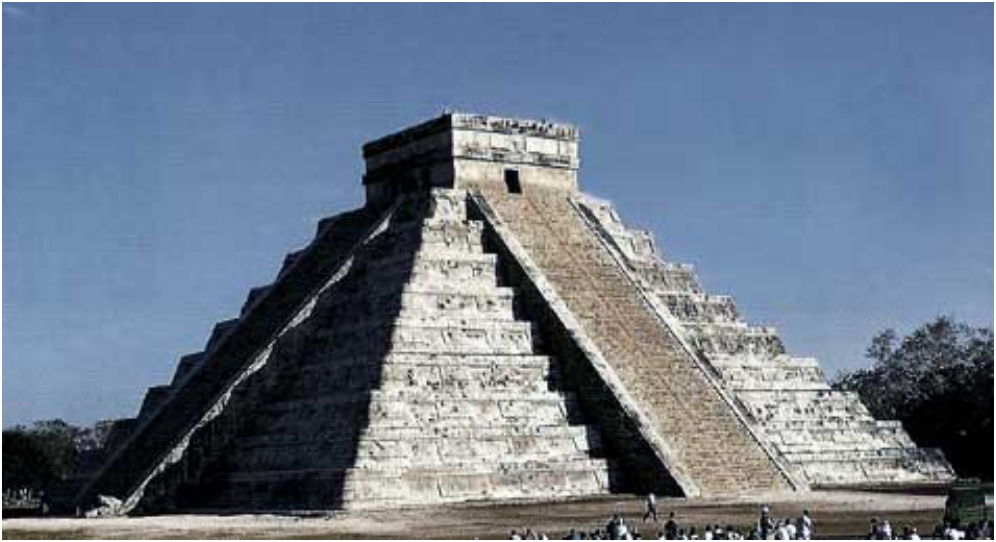}\\
\caption{Pyramid of Kukulkan during an equinox. The pyramid is situated at Chich\'en Itz\'a, Yucat\'an, Mexico.}\label{kukulkan}
\end{figure}

\begin{table}[h!]
\begin{center}
\begin{tabular}{l l l l l l}

Name & $i$ & $C_i$ [day]  & $\mathcal{N}$/$\mathcal{P}_i$ & $\mathcal{Q}_i$ \\\hline
- & 0 & 18 & $\mathcal{N}/13/73/\sum_0^6 C_i$ & 18 \\
Tun & 1 & 360	& $\mathcal{N}/13/73/\sum_0^5 C_i$ & 360\\
Katun & 2 & 7200 & $\mathcal{N}/13/73/\sum_0^4 C_i$ & 7215\\
Baktun & 3 & 144000 & $\mathcal{N}/13/73/\sum_0^3 C_i$ & 144304\\
Pictun & 4 & 2880000 & $\mathcal{N}/13/73/\sum_0^2 C_i$ & 2886428\\
Calabtun & 5 & 57600000 & $\mathcal{N}/13/73/\sum_0^1 C_i$ & 57866020\\
Kinchiltun & 6 & 1152000000 & $\mathcal{N}/13/73/C_0$ & 1215186420\\\hline

\end{tabular}
\caption{Divisibility of the calendar supernumber $\mathcal{N}$ by a polynomial expression of the type $\mathcal{P}_i=13\times73\times (\sum_{n=0}^{6-i}18\times20^n)$. $\mathcal{Q}_i$ is the quotient of the Euclidean division of $\mathcal{N}$ by $\mathcal{P}_i$.}\label{XXX}
\end{center}
\end{table}

Finally, we discuss an important religious site in Mesoamerica, the pyramid of Kukulkan built sometime in the X century CE at Chich\'en Itz\'a (Figure \ref{kukulkan}) where various numbers in the architecture seem to be related to calendar considerations. The pyramid shape may be linked to the Long Count Calendar and the planet canonical cycles which draws a pyramid-like structure (Table \ref{XXX}). The pyramid is constituted of 9 platforms with 4 stairways of 91 steps each leading to the platform temple corresponding to the 3276-day cycle: 3276 = LCM(4,9,91) = 4 $\times$ 819, the coincidence of the 4 directions-colors with the Kawil. The Haab' is represented by the platform temple making the 365$^\mathrm{th}$ step with the 4 $\times$ 91 = 364 steps of the 4 stairways. Each side of the pyramid contains 52 panels corresponding to the Calendar Round: 52 $\times$ 365 = 73 $\times$ 260 = 18980.
The dimensions of the pyramid may be of significance: the width of the top platform is 19.52 m (13 zapal), the height up to the top of the platform temple is 30 m (20 zapal) and the width of the pyramid base is 55.30 m (37 zapal), taking into account the Mayan zapal length measurement such that 1 zapal $\approx$ 1.5 m.\cite{OBrien} The pyramid height and the width of the top platform represents the Tzolk'in (13 $\times$ 20). The stairways divide the 9 platforms of each side of the pyramid into 18 segments which, combined with the pyramid height, represents the 18 Winal of a Tun (18 $\times$ 20). The width of the base 37 and the 4-Kawil 3276 are such that 37 $\times$ 3276 = LCM(148,3276) = 121212 which represents the coincidence of the 4-Kawil and the pentalunex used for solar/lunar eclipse prediction in the Dresden Codex.\cite{bricker} If we include the pentalunex to calculate the calendar supernumber $\mathcal{M}=\mathrm{LCM}(148,\mathcal{N})=37\times\mathcal{N}$ = $36\times\prod_i p_i$, with \{$p_i$, $i$ = 1..11\} = \{2, 3, 5, 7, 13, 19, 29, 37, 59, 73, 89\}, the result of the Euclidean division of $\mathcal{M}$ by $\mathcal{B}=7\times37\times\mathcal{A}=37\times\mathcal{GC}=7\times13\times37\times73\times144000$ can be expressed as a function of the Xultun numbers: 

\begin{eqnarray}\label{astro5}
\mathcal{M}&=&\mathcal{B} \times \mathcal{Q}+ \mathcal{R}\\\nonumber
\mathcal{Q} &=& 21699 \\\nonumber
\mathcal{R} &=& 26810882280\\\nonumber
&=& 2 \times 3^2 \times 7 \times 37\times101 \times56940\\\nonumber
 &= & 7\times18\times37\times\sum_{i=0}^3\mathcal{X}_i.
\end{eqnarray}

During an equinox, the Sun casts a shadow (7 triangles of light and shadow) on the northern stairway representing a serpent snaking down the pyramid (Figure \ref{kukulkan}). The 7 triangles, the 18 segments, the width of the pyramid base (37 zapal) and the Xultun numbers (Table \ref{Xultun}) may be interpreted as a representation of Equ. \ref{astro5}. That gives a connection between the religious sites of Xultun and Chich\'en Itz\'a.

In conclusion, this study presents a unified description of the Mayan Calendar based on naked-eye astronomy. A calendar supernumber $\mathcal{N}$ is calculated by taking the least common multiple of 8 naked-eye astronomy canonical periods describing the planet synodic movements and the apparition of solar/lunar eclipse. This calendar supernumber defines the three major Mayan Calendar cycles: the 4-Kawil cycle, combination of the 4 directions-colors and the 819-day Kawil, the Calendar Round and the Long Count Calendar. The 360-day Tun, the 365-day Haab' and the 3276-day 4-Kawil are issued from the prime factorization of the 8 canonical astronomical input parameters. The 260-day Tzolk'in, the 360-day Tun, the 7200-day Katun and the 144000-day Baktun are obtained from arithmetical calculations on $\mathcal{N}$. The correlation of the three major calendar cycles represents a calendar grand cycle $\mathcal{GC}$. The two Euclidean divisions of $\mathcal{N}$ by $\mathcal{GC}$ and $\mathcal{GC}$/7 = 13 $\times$ 73 $\times$ 144000 show the existence of a Maya Epoch constituted of 5 Maya Eras of 13 Baktun and lead to the Xultun numbers, dating from the early IX century CE and deciphered by Saturno {\it et al.} in 2012 inside a small room of the extensive Mayan ruins of Xultun, Guatemala.\cite{saturno} Modular arithmetic considerations on the calendar supernumber determine the position of the Calendar Round at the Mayan mythical date of creation 0.0.0.0.0 4 Ahau 8 Cumku, reflecting the Mayan cyclical concept of time. The results show a connection between the religious sites of Xultun, Guatemala and Chich\'en Itz\'a, Mexico. This study constitutes a breakthrough towards the understanding of  Mayan ethnomathematics of divination, used to correlate ritual cycles with astronomical events in order to rythm political life and religious practices, embedding Maya civilization in the endless course of time.\\\\
$^*$ e-mail: \texttt{thomas.chanier@gmail.com}

\end{document}